\documentclass[12pt,twoside]{article}

	\usepackage{fancyhdr} 
	\usepackage[margin=1in]{geometry}  
	\usepackage{graphicx}              
	\usepackage{amsmath}               
	\usepackage{amsfonts}              
	\usepackage{amsthm}                
	\usepackage{enumitem}
	\usepackage{amssymb}
	\usepackage{mathtools}
	\usepackage{hyperref}
	\usepackage{listings}
	\usepackage{xcolor}
	\usepackage{booktabs}
	\usepackage{tikz-cd}
	\usepackage{amscd}
	
	\ExplSyntaxOn
	\cs_new_eq:NN \intfromalph \int_from_alph:n
	\ExplSyntaxOff

	\newtheorem{thm}{Theorem}[section]
	\newtheorem{lem}[thm]{Lemma}
	\newtheorem{prop}[thm]{Proposition}
	\newtheorem{cor}[thm]{Corollary}
	\newtheorem{conj}[thm]{Conjecture}
	\theoremstyle{remark}
	\newtheorem{rem}[thm]{Remark}
	
	\theoremstyle{definition}
	\newtheorem{Def}[thm]{Definition}
	\theoremstyle{plain}

	\newtheorem{exa}[thm]{Example}
\begin{document}

		\nocite{*}
		
		\title{Hecke Saturations, Component Groups, and Level Lowering
		}
		
		\author{Mohammad Masih Hamidi\thanks{E-mail: mohammad.hamidi@mail.mcgill.ca}\\ 
			Department of Mathematics \\
			McGill University\footnote{Address: 805 Sherbrooke Street West, Montreal, QC, Canada}\\
		}
		\date{}
		\maketitle
		
\begin{abstract}
	Let $J$ be an abelian variety with split purely toric reduction
	over a non-archimedean local field $F$, equipped with an action
	of a finite flat commutative $\mathbb Z$-algebra $\mathbb T$.
	Under a rational rank-one hypothesis and a local self-duality
	hypothesis, we prove that localized $\mathbb T$-equivariant
	endomorphisms of the rigid uniformizing lattice arise from
	algebraic endomorphisms whenever $\ell\nmid\#\mu(F)$, where
	$\ell$ is the residue characteristic of the localization.
	
	We apply this comparison to optimal $\mathbb T$-quotients,
	extending a component-group cokernel formula of Papikian and
	Rabinoff to higher-dimensional quotients. We also compare the character and component groups of modular Jacobians of level $N$ with those of their $p$-new subvarieties
	for $p\parallel N$, obtaining an alternative proof of Mazur's
	principle when $\ell^2\nmid N$.
	
	For Shimura-curve Jacobians and their new subvarieties, we
	establish component-group surjectivity under suitable local
	hypotheses. We prove a domino principle whereby
	controllability at a single prime of the new level
	implies surjectivity at every prime of that level.
	These results motivate a surjectivity conjecture at
	non-Eisenstein maximal ideals.
\end{abstract}

		\section{Introduction}

Let $F$ be a non-archimedean local field,
$J/F$ be an abelian variety with split purely toric reduction.
Rigid-analytic uniformization realizes $J$ as a quotient of a split torus
by a discrete period lattice $\Lambda$. Every algebraic endomorphism of
$J$ induces an endomorphism of $\Lambda$, but an endomorphism of the
abstract lattice need not respect its embedding into the torus and
therefore need not descend to an algebraic endomorphism of $J$. This
raises a natural question: when is an endomorphism of the uniformizing
lattice algebraic? We study this question when $J$ is equipped with an
action of a finite flat commutative $\mathbb Z$-algebra $\mathbb T$,
locally at maximal ideals of $\mathbb T$.

For an abelian variety $B/F$ with toric reduction, let $\mathcal X_F(B)$
denote the geometric character group of the toric part of the special
fiber of its N\'eron model, and let $\Phi_F(B)$ denote its geometric
component group. We impose the rational rank-one condition
\[
\mathcal X_F(J)\otimes_{\mathbb Z}\mathbb Q
\simeq \mathbb T\otimes_{\mathbb Z}\mathbb Q
\]
as modules over $\mathbb T\otimes_{\mathbb Z}\mathbb Q$. This condition
is satisfied in the Shimura-curve setting of interest here; see
\cite[Lemma 4.13]{helm2007maps}. Rigid uniformization identifies
$\Lambda$ with $\mathcal X_F(J^\vee)$ and supplies an $F^\times$-valued
period pairing between the character groups of $J$ and $J^\vee$. Let
$\mathfrak n$ be a maximal ideal of $\mathbb T$ of residue
characteristic $\ell$. Our first main result assumes local self-duality
at $\mathfrak n$: after a suitable correction of the Hecke action, a
polarization induces a $\mathbb T_{\mathfrak n}$-equivariant isomorphism
between the localized character groups of $J$ and $J^\vee$(Section 2.2). Under this
hypothesis, we obtain the following.

\begin{thm}[Localized endomorphism comparison]
	Assume that $J$ has split purely toric reduction, satisfies the rational
	rank-one property, and is locally self-dual at $\mathfrak n$. If
	$\ell\nmid\#\mu(F)$, then the natural homomorphism
	$
	\operatorname{End}_{\mathbb T}(J)_{\mathfrak n}
	\longrightarrow
	\operatorname{End}_{\mathbb T_{\mathfrak n}}
	(\Lambda_{\mathfrak n})
	$
	is an isomorphism.
\end{thm}

The essential point is that under the rational rank-one hypothesis, a
Hecke-equivariant endomorphism of $\Lambda$ is rationally represented by
an element of the Hecke algebra. Clearing an $\ell$-power denominator
makes the relevant endomorphisms compatible with the rigid period
pairing. Before clearing denominators, the possible discrepancy is an
$\ell$-power root of unity in $F$, which vanishes under the assumption
$\ell\nmid\#\mu(F)$.

This theorem is analogous in spirit to Helm's comparison of
Hecke-equivariant morphisms of Shimura-curve Jacobians with morphisms of
their toric character groups \cite[Theorem 6.1]{helm2007maps}. Our result
applies to a more general class of abelian varieties and uses the rigid
period pairing rather than the controllability and level-raising
arguments of Helm. In particular, it does not require that $\mathfrak n$
be non-Eisenstein. The roots-of-unity condition, however, imposes a
restriction on the field over which the torus splits. In the usual
Shimura-curve setting with trivial nebentypus, the torus splits over an
unramified quadratic extension of $\mathbb Q_p$; our
condition then becomes $\ell\nmid(p^2-1)$. We will use Helm's arithmetic
comparison separately to obtain stronger surjectivity statements in the
Shimura-curve setting.

Our first application concerns component groups of optimal Hecke
quotients. Let $A$ be the factor of $J$ associated with a minimal prime
$I$ of $\mathbb T$ contained in $\mathfrak n$, and let
$\pi:J^\vee\to A^\vee$ be the dual of the optimal inclusion $A\to J$.
Write $R$ for the image of $\mathbb T$ in $\operatorname{End}_F(A)$, let
$\mathfrak m$ be the image of $\mathfrak n$ in $R$, and assume that
$R_{\mathfrak m}$ is a discrete valuation ring $\mathcal O$ with residue
field $k$. Grothendieck's monodromy pairing describes the cokernel of
the map on component groups in terms of the saturation of character
lattices. In Proposition~2.3 we make this description explicit using
the rational Hecke idempotent associated with $A$. If $n$ and $r$ are
its denominator exponents on the $\mathfrak n$-localized $\ell$-adic
Tate module and character group $\mathcal X_F(J^\vee)$, respectively,
we obtain
\[
\#\operatorname{coker}\bigl(
\Phi_F(J^\vee)_{\mathfrak n}
\xrightarrow{\pi_*}\Phi_F(A^\vee)_{\mathfrak m}\bigr)
\leq |k|^{n-r}.
\]
Under local self-duality, this inequality is an equality. In particular, equality of the two
denominator exponents implies surjectivity on localized component
groups. The localized endomorphism comparison in Theorem 1.1 supplies a sufficient
condition for that equality, extending the denominator calculation of
Papikian and Rabinoff \cite{papikian2016optimal} from elliptic quotients
to the higher-dimensional setting considered here.

The second part of the paper considers level lowering. Suppose that
$p\parallel N$, and let $J^{p\text{-new}}=J_0(N)^{p\text{-new}}
\subseteq J_0(N)$ be the $p$-new subvariety. We first compare the
character and component groups of $J_0(N)$ with those of
$J^{p\text{-new}}$. Theorem~4.1 shows that inclusion induces a canonical
Hecke-equivariant isomorphism on geometric character groups at $p$,
and, after localization away from the $p$-old locus, the corresponding
dual quotient induces an isomorphism on component groups. The same
comparison also gives local self-duality of the $p$-new subvariety at
such maximal ideals. This isolates component-group surjectivity as the
input for our level-lowering criterion. More precisely, let $f$ be a
$p$-new normalized eigenform of level $N$, with associated modular
abelian variety $A_f$, and let $\mathfrak m$ be the maximal ideal
determined by its residual Hecke eigensystem. If the residual Galois
representation attached to $\mathfrak m$ is absolutely irreducible and
finite at $p$, Theorem~5.1 shows that surjectivity of the map on
component groups from ${J^{p\text{-new}}}^\vee$ to $A_f^\vee$ forces the
corresponding maximal ideal of the full Hecke algebra to be $p$-old.
Combining this criterion with available character-group
multiplicity-one results gives an alternative approach to Mazur's principle whenever $\ell^2 \nmid N$.

We then turn to the component-group surjectivity problem for optimal
quotients of Jacobians of Shimura curves. Suppose that $N=DM$ is
squarefree, where $D$ is a product of an even number of primes, and
write $J^D(M)$ for the Jacobian of the Shimura curve of discriminant
$D$ and level $M$. For a divisor $Q$ of $M$, let $J_Q^D(M)$ be its
$Q$-new subvariety, and let $A^D(M)^\vee$ be the optimal quotient associated
with a weight-two newform of level $N$. For $D=1$ this includes modular
Jacobians. Ribet and Takahashi \cite{ribet1997parametrizations} proved that for
$p\mid D$ the map on component groups induced by an optimal elliptic
quotient of $J^D(M)$ is surjective locally at non-Eisenstein maximal
ideals when $M$ is not prime. Under the local discrete valuation ring
hypothesis, which is automatic for elliptic quotients, we remove this
restriction for residue characteristics greater than $2$ and extend
the surjectivity result to $Q$-new subvarieties. For $\ell>3$, the
result at primes dividing $D$ also holds with the auxiliary level
structures of \cite{helm2007maps}; the case $\ell=3$ uses the
squarefree $\Gamma_0$-level setting. We also treat residue
characteristic $2$ under a controllability hypothesis.

For primes dividing $Q$, a different phenomenon arises. The character-group comparisons of Helm \cite{helm2007maps}, building on the crucial Theorem 4.1 of Ribet in \cite{ribet1990modular}, relate the character groups
of different Shimura-curve Jacobians at different primes. 
These comparisons allow a local condition at one prime to imply surjectivity
at other primes. We refer to this phenomenon as the domino principle. In this respect, our approach is analogous to that of Ribet--Takahashi
\cite{ribet1997parametrizations}.
 
The following theorem summarizes the principal assertions of
Theorem \ref{domino}. Here controllability is in the sense of
\cite[Definition 6.4]{helm2007maps}.

\begin{thm}[Domino principle and surjectivity at the discriminant]
	Let $J=J_Q^D(M)$ with Hecke algebra $\mathbb T=\mathbb T_Q^D(M)$,
	and let $I$ be the ideal defining the quotient $A^D(M)^\vee=J/IJ$. Let
	$\mathfrak m$ be a non-Eisenstein maximal ideal of $\mathbb T$
	containing $I$, with residue characteristic $\ell$. Suppose that the
	completion of $(\mathbb T/I)_{\mathfrak m}$ is a discrete valuation
	ring. For the induced map
	\[
	\pi_*:\Phi_p(J^\vee)\longrightarrow\Phi_p(A^D(M)^\vee),
	\]
	one has $\operatorname{coker}(\pi_*)_{\mathfrak m}=0$ in each of the
	following cases:
	\begin{enumerate}
		\item $p\mid D$, provided either $\ell>2$ or the inverse image of
		$\mathfrak m$ in the Hecke algebra of $J^D(M)$ is controllable
		at some prime dividing $D$;
		\item $p\mid Q$, provided that $\mathfrak m$ is controllable at
		some prime $q\mid Q$;
		\item $J=J_p^D(M)$ and $p\mid M$, provided that there exists
		a prime $q\mid M$ such that $\mathfrak m$ is not $q$-old.
	\end{enumerate}
	For $\ell>3$, all three cases also hold with auxiliary level
	structure $\Gamma_0(M)\cap\Gamma'$ as in \cite{helm2007maps},
	where the level of $\Gamma'$ is prime to $DM$, for the
	corresponding optimal quotients.
\end{thm}

For odd $\ell$, the condition $\ell\nmid(q-1)$ implies controllability
at $q$ (Remark \ref{rem}). Thus a single prime $q\mid Q$ satisfying this
numerical condition forces surjectivity at every prime dividing $Q$.
In particular, the domino principle applies even when
$\ell\mid(q+1)$, a case excluded by the roots-of-unity condition in Theorem 1.1. The controllability formulation, in this particular case, is more flexible
than the numerical criterion of Theorem 1.1 and also permits certain cases in
residue characteristic $2$. The theorem likewise yields surjectivity
at $p\mid D$ in residue characteristic $2$ when the stated
controllability condition holds. The third case follows by comparing
the component groups of the $p$-new and $pq$-new subvarieties,
reducing surjectivity to the second case of the theorem.

The main mechanism behind the domino principle is the comparison of
character groups supplied by Jacquet--Langlands and Helm's results,
combined with the denominator formula of Proposition~2.3.
Controllability at $q$ gives algebraicity of localized lattice
endomorphisms at that prime, and the character-group comparisons
transfer this algebraicity to the character group at the independently
chosen prime $p$. Thus the arithmetic hypothesis at $q$ propagates
surjectivity to the other relevant primes. 

Motivated by these results and the earlier work of Ribet and Takahashi,
we conjecture (Conjecture \ref{conj}) that, for the optimal quotients of
Shimura-curve Jacobians and their new subvarieties in the more general congruence-level setting of \cite{helm2007maps}, the induced map on component groups at primes dividing the discriminant or the new level is surjective after localization at every non-Eisenstein maximal ideal. The non-Eisenstein restriction is essential in view of the Eisenstein counterexamples in \cite{Hamidi2026Comparison}.

It is natural to ask whether the same conjectural surjectivity holds
for the more general class of abelian varieties with purely toric
reduction considered in Section~2, when the varieties are defined over
$\mathbb Q$. In \cite[Section 5]{papikian2016optimal}, Papikian and
Rabinoff construct optimal elliptic quotients of genus-two Jacobians
for which surjectivity on component groups fails. Their explicit
examples have rational $2$-torsion, so the corresponding mod-$2$
Galois representations are reducible. In Example \ref{exa} we modify their
gluing construction to obtain, for every odd prime $p$, infinitely many
examples of a genus-two Jacobian $J/\mathbb Q$ with two optimal elliptic
quotients $\pi_i:J\to E_i$ ($i=1,2$) satisfying
$
\Phi_p(J)=0$
and $\Phi_p(E_i)\simeq\mathbb Z/2\mathbb Z,
$
where $E_i[2]$ is absolutely irreducible as a $G_{\mathbb Q}$-module.
For the endomorphism order
$
\mathbb T=\{(a,b)\in\mathbb Z\times\mathbb Z:
a\equiv b\pmod 2\},
$
acting on $J$, there is a non-Eisenstein maximal ideal $\mathfrak m$
of residue characteristic $2$ such that
$
\operatorname{coker}((\pi_i)_*)_{\mathfrak m}
\simeq\mathbb F_2~(i=1,2).
$
Moreover, $J$ satisfies the rational rank-one hypothesis of Section~2.
This shows that absolute irreducibility of the residual Galois
representation is not sufficient, by itself, to guarantee surjectivity
on component groups for general abelian varieties with purely toric
reduction. In particular, a nontrivial $2$-primary cokernel need not be
explained by rational $2$-torsion on the optimal quotient. 

The paper is organized as follows. Section~2 develops the monodromy
calculation and the denominator formula for component-group cokernels (Proposition \ref{prop:component-cokernel}).
Section~3 proves the localized endomorphism comparison (Theorem \ref{thm:rigid-endomorphism-comparison}) using
rigid-analytic uniformization. Section~4 compares the character and
component groups of a modular Jacobian and its $p$-new subvariety.
Section~5 establishes the component-group criterion for level lowering,
recovers the stated cases of Mazur's principle, and proves the
surjectivity theorems and domino principle for new subvarieties of
Shimura-curve Jacobians. We conclude with a conjecture on localized
surjectivity and the construction of non-Eisenstein counterexamples
outside our setting.

\section{Character groups and component-group cokernels}
\label{sec:character-cokernel}

The goal of this section is to express the cokernel of the map on
component groups of an optimal Hecke quotient through the denominators of a rational Hecke idempotent. 
Under a local DVR and self-duality hypotheses, we formulate a higher dimensional analog of a formula of Papikian and
Rabinoff \cite[formula 3.7]{papikian2016optimal}.

\subsection{Rational Hecke modules and monodromy}

Let $F$ be a non-archimedean local field and let
$J/F$ be an abelian variety with purely toric reduction. Let $\mathbb T$
be a finite flat commutative $\mathbb Z$-algebra acting faithfully on
$J$. For any abelian variety $B/F$ with purely toric reduction, denote
by $\mathcal X_F(B)$ the geometric character group of the toric part
of the special fiber of its N\'eron model and by $\Phi_F(B)$ its
geometric component group.
We write
$
\mathbb T_{\mathbb Q}=\mathbb T\otimes_{\mathbb Z}\mathbb Q,
~
X_0=\mathcal X_F(J),~
X_0^\vee=\mathcal X_F(J^\vee).
$
As in the introduction, we impose the rational rank-one property

\begin{equation}\label{Rankone}
	X_0\otimes_{\mathbb Z}\mathbb Q\cong\mathbb T_{\mathbb Q}
\end{equation}
as $\mathbb T_{\mathbb Q}$-modules. In the modular and Shimura-curve settings motivating the paper, the
relevant rank-one statement holds (See
\cite[Lemma~4.13]{helm2007maps}).  Let $	\operatorname{End}_{\mathbb T}(J)$ denote the $\mathbb{T}$-equivariant endomorphisms of $J$ and $	\operatorname{End}_{\mathbb T_{\mathbb Q}}
(X_0\otimes\mathbb Q)$ is defined similarly.

\begin{lem}\label{lem:rational-character-rank-one}
	$
	\operatorname{End}_{\mathbb T}(J)\otimes_{\mathbb Z}\mathbb Q
	\;=\;\mathbb T_{\mathbb Q}.
	$
\end{lem}

\begin{proof}
	An algebraic endomorphism of $J$ acts faithfully on $X_0$, as an endomorphism because of the purely toric reduction. Thus restriction gives
	an injection
	$
	\operatorname{End}_{\mathbb T}(J)\otimes\mathbb Q
	\hookrightarrow
	\operatorname{End}_{\mathbb T_{\mathbb Q}}
	(X_0\otimes\mathbb Q).
	$
	By rational rank one assumption, the ring on the right is $\mathbb T_{\mathbb Q}$.
	The reverse inclusion follows from the faithful Hecke action.
\end{proof}

Let
$
\pi:J^\vee\twoheadrightarrow A^\vee
$
be a $\mathbb{T}$-equivariant optimal quotient, where $A^\vee$ is
a simple abelian variety, and let
$
i=\pi^\vee:A\hookrightarrow J
$
be the dual inclusion. Put
$
I=\ker\bigl(\mathbb T\longrightarrow
\operatorname{End}(A)\bigr).
$
We assume that the rational projector onto the factor $i(A)$
belongs to $\mathbb T_{\mathbb Q}$, and denote it by $e$.

Put
$
Y_0=\mathcal X_F(A),~
Y_0^\vee=\mathcal X_F(A^\vee),
~
M^*=\operatorname{Hom}_{\mathbb Z}(M,\mathbb Z).
$
Grothendieck's monodromy pairing gives functorial exact sequences
\[
0\longrightarrow X_0\longrightarrow (X_0^\vee)^*
\longrightarrow\Phi_F(J^\vee)\longrightarrow0,
\]
\[
0\longrightarrow Y_0\longrightarrow (Y_0^\vee)^*
\longrightarrow\Phi_F(A^\vee)\longrightarrow0.
\]
The map $\pi$ induces an injection
$
\pi^*:Y_0^\vee\hookrightarrow X_0^\vee,
$
while the closed immersion $i$ induces a surjection
$
i^*:X_0\twoheadrightarrow Y_0.
$
Define the saturation of the image of $\pi^*$ by
$
\overline{Y_0^\vee}
=\bigl(\pi^*Y_0^\vee\otimes_{\mathbb Z}\mathbb Q\bigr)
\cap X_0^\vee
\subset X_0^\vee\otimes\mathbb Q.
$
The commutative monodromy diagram and the snake lemma identify the
cokernel of $\pi_*$ with the cokernel of the middle vertical map. Hence
\begin{equation}\label{eq:general-saturation-cokernel}
	\operatorname{Coker}\bigl(
	\Phi_F(J^\vee)\xrightarrow{\pi_*}\Phi_F(A^\vee)\bigr)
	\cong
	\operatorname{Hom}_{\mathbb Z}
	\left(\overline{Y_0^\vee}/\pi^*Y_0^\vee,
	\mathbb Q/\mathbb Z\right).
\end{equation}
In particular,
\begin{equation}\label{eq:general-saturation-index}
	\#\operatorname{Coker}(\pi_*)
	=\bigl[\overline{Y_0^\vee}:\pi^*Y_0^\vee\bigr].
\end{equation}
This is the general saturation formula of
\cite[Proposition~2.12]{papikian2004component}.
For the Hecke factor cut out by $I$, the saturated lattice is
$
\overline{Y_0^\vee}=X_0^\vee[I]=\{x\in X_0^\vee:Ix=0\}.
$

\subsection{Local self-duality and coefficient orders}

Let $R$ be the image of $\mathbb T$ in $\operatorname{End}_F(A)$.
Since $A$ is simple, $\operatorname{End}_F^0(A)$ is a division
algebra. Hence $R$ is an order in a number field
$K=\operatorname{Frac}(R)$.
 Let $\mathfrak m\subset R$
be a maximal ideal of residue characteristic $\ell$, and write
$\mathfrak n\subset\mathbb T$ for its inverse image. A subscript
$\mathfrak n$ refers to localization of $\mathbb T$-modules, whereas a
subscript $\mathfrak m$ refers to localization of $R$-modules.
Fix a polarization
$
H_1:J\longrightarrow J^\vee.
$

\begin{Def}\label{def:local-self-duality}
	We say that $J$ is \emph{locally self-dual at $\mathfrak n$} if there
	exists an automorphism $\omega:J\to J$ such that the map
	$
	H=H_1\circ\omega:J\longrightarrow J^\vee
	$
	is $\mathbb T$-equivariant and induces an isomorphism
	$
	T_\ell(H)_{\mathfrak n}:T_\ell(J)_{\mathfrak n}
	\xrightarrow{\sim}T_\ell(J^\vee)_{\mathfrak n}
	$
	of $\mathbb T_{\mathfrak n}$-modules.
\end{Def}

The correction $\omega$ permits us to account for the adjoint Hecke
action without requiring the original polarization to be
$\mathbb T$-equivariant. In modular examples an Atkin--Lehner involution facilitates this.
Next, assume that the
coefficient order is maximal at $\mathfrak m$:
\begin{equation}\label{eq:coefficient-DVR}
	R_{\mathfrak m}=(\mathcal O_K)_{\mathfrak P}
	=:\mathcal O,
\end{equation}
where $\mathcal O_K$ is the full ring of integers of $K$ and
$\mathfrak P$ lies above $\mathfrak m$. Let $\mathfrak p$ be the
maximal ideal of $\mathcal O$, let $k=\mathcal O/\mathfrak p$. Normalize
$\operatorname{ord}_{\mathfrak p}$ so that a uniformizer has
valuation one. Write
$
X=(X_0)_{\mathfrak n},~ X^\vee=(X_0^\vee)_{\mathfrak n},
~ Y=(Y_0)_{\mathfrak m},~
Y^\vee=(Y_0^\vee)_{\mathfrak m}.
$
The rational rank-one property identifies $(X_0\otimes\mathbb Q)[I]$ with the coefficient field $K$. Since $Y$ and
$Y^\vee$ are finitely generated torsion-free rank-one modules over
the DVR $\mathcal O$, they are free:
\begin{equation}\label{eq:rank-one-character-lattices}
	Y\cong\mathcal O,~ Y^\vee\cong\mathcal O.
\end{equation}
The monodromy pairing and the rational rank-one property also give
\begin{equation}\label{eq:I-lattice-rank-one}
	X^\vee[I]\cong X[I]\cong\mathcal O.
\end{equation}
Let $e\in\mathbb T_{\mathbb Q}$ denote the rational idempotent
corresponding to the factor $A$. On the $\mathfrak n$-localized
character lattice define the denominator exponent
\[
r=\min\left\{
\operatorname{ord}_{\mathfrak p}(a):0\ne a\in\mathcal O,
\ (ae)(X^\vee)\subseteq X^\vee\right\}.
\]
Likewise, for the $\mathfrak n$-part of the $\ell$-adic Tate module,
put
\[
n=\min\left\{
\operatorname{ord}_{\mathfrak p}(a):0\ne a\in\mathcal O,
\ (ae)T_\ell(J)_{\mathfrak n}
\subseteq T_\ell(J)_{\mathfrak n}\right\}.
\]
One checks that $r$ and $n$ are well-defined. 

\subsection{The higher-dimensional denominator formula}
In Tate-module
calculations below, localized coefficient and character modules are
understood after $\mathfrak p$-adic completion.
\begin{prop}\label{prop:component-cokernel}
	Assume the rational rank-one hypothesis \ref{Rankone}
	and the maximality condition \eqref{eq:coefficient-DVR}.
	With $n$ and $r$ defined in previous section, the following hold.
	\begin{enumerate}
		\item[(i)] Without a local self-duality assumption,
		\[
		\#\operatorname{Coker}\left(
		\Phi_F(J^\vee)_{\mathfrak n}
		\xrightarrow{\pi_*}\Phi_F(A^\vee)_{\mathfrak m}
		\right)\le |k|^{\,n-r};
		\]
		\item[(ii)] If $J$ is additionally locally self-dual at $\mathfrak n$, and $\omega$ restricts to a $\mathbb{T}$-equivariant automorphism of $A$, 
		then equality holds:
		\[
		\#\operatorname{Coker}\left(
		\Phi_F(J^\vee)_{\mathfrak n}
		\xrightarrow{\pi_*}\Phi_F(A^\vee)_{\mathfrak m}
		\right)=|k|^{\,n-r}.
		\]
	\end{enumerate}
\end{prop}

In the classical modular-Jacobian setting, as mentioned after
Definition~\ref{def:local-self-duality}, the correction $\omega$ can
be taken to be an Atkin--Lehner involution. For modular abelian
varieties attached to newforms with trivial nebentypus, its
restriction acts as $\pm 1$. Hence the additional hypothesis
on $\omega$ in part~(ii) of Proposition~\ref{prop:component-cokernel}
is automatically satisfied.
\begin{proof}[Proof of Proposition~\ref{prop:component-cokernel}]
	By \eqref{eq:rank-one-character-lattices} and
	\eqref{eq:I-lattice-rank-one}, choose generators
	$t\in Y^\vee$ and $s\in X^\vee[I]$. There is a nonzero
	$c\in\mathcal O$ such that
	$
	\pi^*(t)=cs.
	$
	The saturation formula \eqref{eq:general-saturation-cokernel}
	therefore gives
	\begin{equation}\label{eq:cokernel-character-index}
		\#\operatorname{Coker}(\pi_*)
		=\#(\mathcal O/c\mathcal O)
		=|k|^{\operatorname{ord}_{\mathfrak p}(c)}.
	\end{equation}
	
	\emph{(i) The general upper bound.}
	Choose $b\in\mathcal O$ with
	$\operatorname{ord}_{\mathfrak p}(b)=n$ such that
	$
	(be)T_\ell(J)_{\mathfrak n}\subseteq T_\ell(J)_{\mathfrak n}.
	$
	Rationally, $e$ is the projector onto the factor
	$i(A)\subseteq J$. Since $i=\pi^\vee$ is a closed immersion,
	the inclusion
	$T_\ell(A)_{\mathfrak m}\hookrightarrow T_\ell(J)_{\mathfrak n}$
	has saturated image. Consequently,
	$
	(be)T_\ell(J)_{\mathfrak n}
	\subseteq i\bigl(T_\ell(A)_{\mathfrak m}\bigr).
	$
	Thus $be$ factors integrally through $T_\ell(A)_{\mathfrak m}$.
	By functoriality of the toric Tate-module filtrations, passing
	to their respective character-lattice quotients gives
	\begin{equation}\label{eq:general-denominator-inclusion}
		(be)X^\vee\subseteq\pi^*Y^\vee.
	\end{equation}
	
	Because $e$ acts as the identity on $X^\vee[I]=\mathcal O s$,
	the definition of $r$ and the DVR hypothesis imply
	$
	e(X^\vee)=\mathfrak p^{-r}\mathcal O s,
	~
	(be)X^\vee=\mathfrak p^{n-r}\mathcal O s.
	$
	Combining this with \eqref{eq:general-denominator-inclusion}
	and $\pi^*Y^\vee=c\mathcal O s$ yields
	$
	\mathfrak p^{n-r}\mathcal O s
	\subseteq c\mathcal O s,
	~
	\operatorname{ord}_{\mathfrak p}(c)\le n-r.
	$
	Part~(i) follows from \eqref{eq:cokernel-character-index}.
	
	\emph{(ii) Equality under local self-duality.}
	Now assume $J$ is locally self-dual at $\mathfrak n$, and that
	$\omega$ restricts to a $\mathbb{T}$-equivariant automorphism
	$u$ of $A$. Write $i=\pi^\vee$, and form
	\[
	q=\pi\circ H:J\longrightarrow A^\vee,
	~
	\alpha=i^\vee\circ H_1\circ i,
	~
	\alpha'=q\circ i=\alpha\circ u.
	\]
	Since $\alpha'$ and $u$ are $\mathbb T$-equivariant, so is
	$\alpha$. Moreover, $\alpha$ is a polarization.
	With $t$ and $s$ fixed above, choose a $\mathcal{O}$-generator $t'\in Y$.
	There are nonzero $a,a_0\in\mathcal O$ such that
	\begin{equation}\label{eq:alpha-c-character}
		\pi^*(t)=cs,\qquad
		(\alpha')^*(t)=at',\qquad \alpha^*(t)=a_0t'.
	\end{equation}
	Since $u$ acts on $Y$ by a unit of $\mathcal O$, the
	elements $a$ and $a_0$ differ by a unit.
	The injection $\pi^*:Y^\vee\hookrightarrow X^\vee[I]$ is the first
	source of the index we want to calculate.
	
	\emph{The Tate-module denominator.}
	The toric Tate-module
	filtrations for $A$ and $A^\vee$ take the form
	\[
	0\longrightarrow
	\operatorname{Hom}_{{\mathcal O}}
	(Y,
	{\mathcal O}(1))
	\longrightarrow T_\ell(A)_{\mathfrak m}
	\longrightarrow Y^\vee
	\longrightarrow0,
	\]
	\[
	0\longrightarrow
	\operatorname{Hom}_{{\mathcal O}}
	(Y^\vee,
	{\mathcal O}(1))
	\longrightarrow T_\ell(A^\vee)_{\mathfrak m}
	\longrightarrow Y
	\longrightarrow0.
	\]
	If $\ell\ne\operatorname{char}(F)$, then
	$T_\ell(A)_{\mathfrak m}$ is free of rank two over
	${\mathcal O}$. With bases adapted to these filtrations and
	the Weil pairing, the alternating form induced by the polarization
	$\alpha$ has matrix
	$
	\begin{pmatrix}0&a_0\\-a_0&0\end{pmatrix}.$
	Since $\alpha'=\alpha\circ u$ and $u$ is an automorphism, the
	cokernels of $T_\ell(\alpha')_{\mathfrak m}$ and
	$T_\ell(\alpha)_{\mathfrak m}$ are isomorphic.
	If $\ell=\operatorname{char}(F)$, then
	$\mathbb Z_\ell(1)=0$ on geometric points, so the filtrations
	give $T_\ell(A)_{\mathfrak m}\cong Y^\vee\cong\mathcal O$ and
	$T_\ell(A^\vee)_{\mathfrak m}\cong Y\cong\mathcal O$.
	Under these identifications $T_\ell(\alpha')_{\mathfrak m}$
	is multiplication by $a$. Thus,
	\begin{equation}\label{eq:tate-alpha-exponent}
		\operatorname{Coker}(T_\ell(\alpha')_{\mathfrak m})
		\cong
		\begin{cases}
			({\mathcal O}/a{\mathcal O})^2,
			&\ell\ne\operatorname{char}(F),\\
			{\mathcal O}/a{\mathcal O},
			&\ell=\operatorname{char}(F).
		\end{cases}
	\end{equation}
	The local self-duality assumption implies that $H$ induces an
	isomorphism on the $\mathfrak n$-local Tate modules. Since $\pi$ is
	optimal, the induced map
	$
	q:T_\ell(J)_{\mathfrak n}\twoheadrightarrow
	T_\ell(A^\vee)_{\mathfrak m}
	$
	is surjective, and the image of
	$i:T_\ell(A)_{\mathfrak m}\hookrightarrow T_\ell(J)_{\mathfrak n}$
	is saturated. Since $q$ and $i$ are Hecke-equivariant and $q\circ i=\alpha'$, we have
	$
	e=i\circ(\alpha')^{-1}\circ q.
	$
	It follows that, for $b\in\mathcal O$, the operator $be$ preserves
	$T_\ell(J)_{\mathfrak n}$ exactly when
	$
	b(\alpha')^{-1}T_\ell(A^\vee)_{\mathfrak m}
	\subseteq T_\ell(A)_{\mathfrak m}.
	$
	The elementary divisors in \eqref{eq:tate-alpha-exponent} therefore
	give
	\begin{equation}\label{eq:denominator-n}
		n=\operatorname{ord}_{\mathfrak p}(a).
	\end{equation}
	
	\emph{The character-group denominator.}
	Since $i:A\hookrightarrow J$ is a closed immersion, the map
	$i^*:X\twoheadrightarrow Y$ is surjective. 	By functoriality of the toric Tate-module filtration,
	the local self-duality assumption also implies that the induced map
	$H^*:X^\vee \rightarrow X$ is an isomorphism, and thus
	$
	q_X=i^*\circ H^*:X^\vee\twoheadrightarrow Y
	$
	is surjective as well. The identity
	$(\alpha')^*=q_X\circ\pi^*$ and \eqref{eq:alpha-c-character}
	give
	$
	q_X(s)=(a/c)t',
	~
	\operatorname{ord}_{\mathfrak p}(a)
	\ge\operatorname{ord}_{\mathfrak p}(c).
	$
	On $X^\vee\otimes\mathbb Q$, the idempotent has the expression
	$
	e=\pi^*\circ((\alpha')^*)^{-1}\circ q_X,
	~
	e(X^\vee)=(c/a)\mathcal O s.
	$
	Since $s$ generates the saturated rank-one $\mathcal O$-lattice
	$X^\vee[I]$, we conclude that
	\begin{equation}\label{eq:denominator-r}
		r=\operatorname{ord}_{\mathfrak p}(a)
		-\operatorname{ord}_{\mathfrak p}(c).
	\end{equation}
	Combining \eqref{eq:denominator-n} and \eqref{eq:denominator-r}
	shows that
	$\operatorname{ord}_{\mathfrak p}(c)=n-r$.
	Therefore, by \eqref{eq:cokernel-character-index},
	\[
	\#\operatorname{Coker}(\pi_*)
	=|k|^{\,n-r},
	\]
	as required. 
	\end{proof}

		\section{Rigid analytic endomorphisms and local surjectivity}
\label{sec:rigid-endomorphisms}

The purpose of this section is to study the algebraicity of Hecke-equivariant
endomorphisms of the period lattice of an abelian variety with split toric
reduction. We therefore begin with a maximal ideal of the
\emph{ambient} Hecke algebra.
Let $F$ be a non-archimedean local field, and let
$J/F$ be an abelian variety with split purely toric reduction, equipped
with a faithful action of a finite flat commutative $\mathbb Z$-algebra
$\mathbb T$. As in Section~2, assume the rational rank-one property
$
\mathcal X_F(J)\otimes_{\mathbb Z}\mathbb Q
\cong \mathbb T\otimes_{\mathbb Z}\mathbb Q.
$
Fix a maximal ideal $\mathfrak n\subset\mathbb T$ of residue
characteristic $\ell$.

\subsection{Split toric uniformization and algebraicity}

Write the rigid-analytic uniformizations of $J$ and $J^\vee$ as
\begin{equation}\label{uniformization}
	0\longrightarrow\Lambda_0
	\overset{j}{\longrightarrow}G^{\mathrm{an}}
	\longrightarrow J^{\mathrm{an}}\longrightarrow0,
\end{equation}
\[
0\longrightarrow\Lambda_0^\vee
\overset{j^\vee}{\longrightarrow}(G^\vee)^{\mathrm{an}}
\longrightarrow(J^\vee)^{\mathrm{an}}\longrightarrow0,
\]
where $G$ and $G^\vee$ are split tori over $F$. The uniformizations are related by the $F^\times$-valued period pairing
\[
q:\Lambda_0\times\Lambda_0^\vee\longrightarrow F^\times,
\qquad
j^\vee(y)(x)=q(x,y).
\]
The character group of $G^\vee$ is $\Lambda_0$, so an integral
endomorphism of $\Lambda_0$ induces an algebraic endomorphism of
$G^\vee$.
Write
$
\Lambda_{\mathfrak n}=(\Lambda_0)_{\mathfrak n},
~
\Lambda^\vee_{\mathfrak n}=(\Lambda_0^\vee)_{\mathfrak n}.
$
By the theory of rigid-analytic uniformization (see \cite[Exp.~IX, \S7, \S14]{grothendieck1972modeles}), there are canonical and functorial identifications $\Lambda_0\cong\mathcal X_F(J^\vee),
~
\Lambda_0^\vee\cong\mathcal X_F(J)
$. Therefore, assuming that $J$ is locally self-dual at $\mathfrak n$, in the sense of
Definition~2.2, there is a global $\mathbb T$-equivariant map
$
H_0:\Lambda_0\longrightarrow\Lambda_0^\vee
$
which induces a $\mathbb T_{\mathfrak n}$-linear isomorphism
$
H:\Lambda_{\mathfrak n}
\overset{\sim}{\longrightarrow}
\Lambda^\vee_{\mathfrak n}.
$
The Hecke actions on the two lattices are adjoint with respect to $q$: for $t\in\mathbb T$ one has
$
q(tx,y)=q(x,ty)
~(x\in\Lambda_0,\ y\in\Lambda_0^\vee),
$
with the dual Hecke action understood on the second variable.
 Let $\mu(F)$ denote the subgroup of roots of unity in $F^{\times}$.
\begin{thm}[Localized endomorphism comparison]
	\label{thm:rigid-endomorphism-comparison}
	 If
	$
	\ell\nmid\#\mu(F),
	$
	then the natural homomorphism
	$
	\operatorname{End}_{\mathbb T}(J)_{\mathfrak n}
	\longrightarrow
	\operatorname{End}_{\mathbb T_{\mathfrak n}}
	(\Lambda_{\mathfrak n})
	$
	is an isomorphism. 
\end{thm}

\begin{rem}[Connection with a result of Helm]
	Theorem~\ref{thm:rigid-endomorphism-comparison} is analogous in spirit to
	\cite[Theorem~6.1]{helm2007maps}, which compares Hecke-equivariant morphisms of
	Shimura-curve Jacobians with morphisms of their toric character groups.
	Our result gives a comparison of localized Hecke-equivariant endomorphism
	rings for a more general class of abelian varieties, assuming the rational
	rank-one property and localized self-duality. Its proof uses rigid analytic
	uniformization and the period pairing.
	Helm introduces a local condition called \emph{controllability}, which
	ensures that the relevant localized character group is free of rank one
	over the Hecke algebra; see \cite[Definition~6.4 and Lemma~6.5]{helm2007maps}.
	In the usual toric setting, the numerical condition $\ell\nmid p-1$
	is sufficient for controllability. Helm proves the desired comparison
	under controllability in $loc. cit.$, and then removes
	this assumption in Theorem~6.1 by a level-raising argument. Our approach does not use controllability, but requires
	$\ell\nmid\#\mu(F)$, where $F$ is the field over which the torus
	splits. In the standard Shimura-curve setting with trivial nebentypus, Frobenius acts on the toric character group through an
	Atkin--Lehner involution; hence the torus splits over the unramified
	quadratic extension $F=\mathbb Q_{p^2}$ of $\mathbb Q_p$.
	Our condition then becomes $\ell\nmid p^2-1$.
	This is stronger than $\ell\nmid p-1$: in particular, it also excludes
	$\ell\mid p+1$. Unlike Helm's Theorem~6.1, which is established away from
	an exceptional set containing the Eisenstein maximal ideals,
	our comparison imposes no non-Eisenstein condition on
	$\mathfrak n$. In particular, it remains valid at Eisenstein
	maximal ideals.
\end{rem}
\begin{proof}
	Injectivity follows from the faithfulness of the action of $\operatorname{End}(J)$ on $\Lambda_0$. Because $\Lambda_0$ is a finitely presented $\mathbb T$-module,
	localization commutes with the endomorphism module, i.e.
	$
	\operatorname{End}_{\mathbb T_{\mathfrak n}}
	(\Lambda_{\mathfrak n})
	\cong\operatorname{End}_{\mathbb T}(\Lambda_0)_{\mathfrak n}.
	$
	It therefore suffices to show that every
	$\psi\in\operatorname{End}_{\mathbb T}(\Lambda_0)$ becomes algebraic
	after localization at $\mathfrak n$.
	Rationally, the rank-one hypothesis gives
	$
	\operatorname{End}_{\mathbb T}(J)\otimes\mathbb Q
	=\mathbb T\otimes\mathbb Q
	=\operatorname{End}_{\mathbb T\otimes\mathbb Q}
	(\Lambda_0\otimes\mathbb Q).
	$ Thus, there exist a positive integer $d$ and
	$t\in\mathbb T$ such that $d\psi=t$ on $\Lambda_0$.
	Write $d=u\ell^a$ with $\gcd(u,\ell)=1$, and put
	$\eta=u\psi$. Since $u$ is a unit in $\mathbb T_{\mathfrak n}$,
	it is enough to prove that $\eta$ is algebraic after localization.
	We have $\ell^a\eta=t\in\mathbb T$.
	For $x,z\in\Lambda_0$, $\mathbb{T}$-adjointness of the period pairing and
	$\mathbb T$-equivariance of $H$ give
	\[
	\begin{aligned}
		q(\eta x,Hz)^{\ell^a}
		&=q(tx,Hz)\\
		&=q(x,tHz)\\
		&=q(x,H(tz))
		=q(x,H(\eta z))^{\ell^a}.
	\end{aligned}
	\]
	Hence, $q(\eta x,Hz)/q(x,H(\eta z))$ lies in $\mu_{\ell^a}(F)$, which is trivial by hypothesis. Consequently
	$
	q(\eta x,Hz)=q(x,H(\eta z))
	.
	$
	By definition, $G^\vee=\text{Hom}(\Lambda_0,\text{G}_m)$. Let $\eta^G$ be the endomorphism of $G^\vee$ induced by $\eta$. The preceding identity says exactly
	that
	$
	\eta^G\bigl(j^\vee(Hz)\bigr)
	=j^\vee\bigl(H(\eta z)\bigr),
	$
	for every $z \in \Lambda_0$. Set $M=\Lambda_0^\vee$ and let
	$
	\delta:
	\operatorname{End}_{\mathbb T}(\Lambda_0)
	\oplus\operatorname{End}_{\mathbb T}(M)
	\longrightarrow
	\operatorname{Hom}_{\mathbb T}(M,G^\vee(F))
	$
	be the homomorphism defined by
	$
	\delta(v,w)=v^G\circ j^\vee-j^\vee\circ w
	$
	($v^G$ is the induced map from $v$ on $G^\vee$).
	The kernel consists of pairs of maps that descend
	to $\mathbb T$-equivariant rigid-analytic endomorphisms
	of $(J^\vee)^{\mathrm{an}}$. Since $J^\vee$ is projective,
	these endomorphisms are algebraic by rigid-analytic GAGA. Since localization is
	exact, the localized kernel is the kernel of $\delta_{\mathfrak n}$.
	Because $H$ is an isomorphism, the formula
	$
	\eta_M=H\eta_{\mathfrak n}H^{-1}
	$
	defines an element of
	$\operatorname{End}_{\mathbb T_{\mathfrak n}}(M_{\mathfrak n})$.
	The identity $\eta^G\bigl(j^\vee(Hz)\bigr)
	=j^\vee\bigl(H(\eta z)\bigr)$ and the equality
	$H(\Lambda_0)_{\mathfrak n}=M_{\mathfrak n}$ show that
	$
	\delta_{\mathfrak n}(\eta,\eta_M)=0.
	$
    Hence, $(\eta,\eta_M)$ lies in the localization of kernel of $\delta$.
	It follows that $\eta$ is induced, after localization, by an element
	of $\operatorname{End}_{\mathbb T}(J^\vee)_{\mathfrak n}$.
	Dualizing gives the desired element of
	$\operatorname{End}_{\mathbb T}(J)_{\mathfrak n}$. The construction above is canonical, hence gives an inverse to the inclusion $
	\operatorname{End}_{\mathbb T}(J)_{\mathfrak n}
	\longrightarrow
	\operatorname{End}_{\mathbb T_{\mathfrak n}}
	(\Lambda_{\mathfrak n})
	$.
\end{proof}

\subsection{Idempotent denominators and component groups}
\label{subsec:denominator-application}
Here, we see an application of previous sections. Let us return to the optimal Hecke quotient from Section~2,
$
\pi:J^\vee\twoheadrightarrow A^\vee,
~ i=\pi^\vee:A\hookrightarrow J.
$
Let $R=\mathbb T_A$ be the coefficient order of $A$, let
$\mathfrak m\subset R$ be a maximal ideal of residue characteristic
$\ell$, and let $\mathfrak n\subset\mathbb T$ be its inverse image. Let $e\in\mathbb T\otimes\mathbb Q$ be the idempotent attached to $A$.
Denote by $n$ and $r$ the denominator exponents of $e$ on the
$\mathfrak n$-localized $\ell$-adic Tate module of $J$ and its
$\mathfrak n$-localized character group
$\mathcal X_F(J^\vee)$, respectively, as in Section~2.
Proposition~\ref{prop:component-cokernel} gives
$
\#\operatorname{Coker}\bigl(
\Phi_F(J^\vee)_{\mathfrak n}
\longrightarrow\Phi_F(A^\vee)_{\mathfrak m}\bigr)
=|k|^{\,n-r}.
$

\begin{cor}[Surjectivity via the endomorphism comparison]\label{cor:pal-type-surjectivity}
	Assume the rational rank-one hypothesis \ref{Rankone} and the
	maximality condition \ref{eq:coefficient-DVR}. Suppose, moreover, that $J$
	has split purely toric reduction over $F$, is locally
	self-dual at $\mathfrak n$, and
	$\ell\nmid\#\mu(F)$.
	Let $n$ and $r$ be the denominator exponents of
	the rational Hecke idempotent $e$ on
	$T_\ell(J)_{\mathfrak n}$ and
	$\mathcal X_F(J^\vee)_{\mathfrak n}$, respectively,
	as defined in Section~2.2.
	Then $n=r$, and the map
	\[
	\Phi_F(J^\vee)_{\mathfrak n}
	\longrightarrow
	\Phi_F(A^\vee)_{\mathfrak m}
	\]
	is surjective.
\end{cor}

To prove the corollary, we need the following lemma.

\begin{lem}\label{lem:tate-character-integrality}
	Let $u\in\mathbb T_{\mathfrak n}[1/\ell]$, viewed in the rational
	Hecke action on $J$. If $u$ preserves the localized $\ell$-adic Tate
	module $T_\ell(J)_{\mathfrak n}$, then it preserves
	$\mathcal X_F(J^\vee)_{\mathfrak n}$.
\end{lem}

    \begin{proof}
    	Let $G$ be the uniformizing torus of $J$, and let $\Lambda_0$ be its period lattice. Applying the Tate-module functor to the exact sequence \ref{uniformization}, we obtain a functorial exact sequence
    	\[
    	0\longrightarrow T_\ell(G)\longrightarrow T_\ell(J)
    	\longrightarrow\Lambda_0\otimes_{\mathbb Z}\mathbb Z_\ell
    	\longrightarrow0.
    	\]
    	Using the canonical identification $\Lambda_0\cong\mathcal X_F(J^\vee)$ and localizing at $\mathfrak n$, we obtain
    	\[
    	0\longrightarrow T_\ell(G)_{\mathfrak n}
    	\longrightarrow T_\ell(J)_{\mathfrak n}
    	\longrightarrow\mathcal X_F(J^\vee)_{\mathfrak n}
    	\longrightarrow0.
    	\]
    	By functoriality, $u$ preserves the rational subspace $V_\ell(G)_{\mathfrak n}$. Since the quotient is $\mathbb Z_\ell$-torsion-free,
    	$
    	T_\ell(G)_{\mathfrak n}
    	=T_\ell(J)_{\mathfrak n}\cap V_\ell(G)_{\mathfrak n}.
    	$
    	Thus, if $u$ preserves $T_\ell(J)_{\mathfrak n}$, it also preserves $T_\ell(G)_{\mathfrak n}$ and hence induces an endomorphism of the quotient $\mathcal X_F(J^\vee)_{\mathfrak n}$.
    \end{proof}

\begin{proof}[Proof of Corollary \ref{cor:pal-type-surjectivity}]
	For $c\in\mathcal O$, the rational operator $ce$ is viewed in
	$\mathbb T_{\mathfrak n}[1/\ell]$ through a lift of $c$ from $\mathcal{O}$; its action on the $A$-factor is independent of the
	chosen lift. By Lemma~\ref{lem:tate-character-integrality},
	$
	ce\in\operatorname{End}_{\mathbb T_{\mathfrak n}}
	(T_\ell(J)_{\mathfrak n})
	$ implies that
	$ce\in\operatorname{End}_{\mathbb T_{\mathfrak n}}
	(\mathcal X_F(J^\vee)_{\mathfrak n}).
	$
	Conversely, if $ce$ preserves the latter character group, Theorem~
	\ref{thm:rigid-endomorphism-comparison} identifies it with an element
	of $\operatorname{End}_{\mathbb T}(J)_{\mathfrak n}$, so it acts
	integrally on $T_\ell(J)_{\mathfrak n}$. The two denominator ideals
	therefore coincide, giving $n=r$. The claim now follows from
	Proposition~\ref{prop:component-cokernel}.
\end{proof}

		\section{The character and component groups of the $p$-new subvariety}
\label{sec:p-new-comparison}

Let $p\parallel N$ and write
$
J=J_0(N),~ J^{p\text{-new}}=J_0(N)^{p\text{-new}}\subseteq J.
$
The latter variety is the $p$-new subvariety. We regard these abelian varieties as defined over $\mathbb Q_p$ by base
change from $\mathbb{Q}$. The variety $J^{p\text{-new}}$ has purely toric reduction at $p$.
For an abelian variety $B/\mathbb Q_p$ with semistable reduction, write
$\mathcal X_p(B)$ for the geometric character group of the toric part of
its N\'eron model and $\Phi_p(B)$ for its component group.
Let $\mathbb T(N)$ be the integral Hecke algebra acting on $J$ and put
$
\mathbb T^{p\text{-new}}(N)
=\operatorname{Im}\bigl(\mathbb T(N)\longrightarrow
\operatorname{End}_{\mathbb Q}(J^{p\text{-new}})\bigr).
$
Write $i:J^{p\text{-new}}\hookrightarrow J$ for the
inclusion. 

\begin{thm}\label{thm:pnew-component}
	The inclusion $i$ induces a canonical $\mathbb T(N)$-equivariant
	isomorphism
	$
	i^*:\mathcal X_p(J)\xrightarrow{\ \sim\ }
	\mathcal X_p(J^{p\text{-new}}).
	$
	Let $\mathfrak m$ be a maximal ideal of
	$\mathbb T^{p\text{-new}}(N)$, and let $\mathfrak n$ denote its preimage
	in $\mathbb T(N)$. If $\mathfrak n$ is not $p$-old, then
	$J^{p\text{-new}}$ is locally self-dual at $\mathfrak m$, and
	$
	\Phi_p(J^\vee)_{\mathfrak n}
	\xrightarrow{\ \sim\ }
	\Phi_p((J^{p\text{-new}})^\vee)_{\mathfrak m}
	$
	via $i^\vee$.
\end{thm}
Let
$
I=\ker\bigl(\mathbb T(N)\longrightarrow
\mathbb T^{p\text{-new}}(N)\bigr).
$
Let $\phi:J\xrightarrow{\sim}J^\vee$ be the canonical principal
polarization. As usual, use the Atkin--Lehner involution to correct the
adjoint Hecke action and denote the resulting $\mathbb T(N)$-equivariant
isomorphism by
$
\lambda:J\xrightarrow{\ \sim\ }J^\vee
$
($i.e.~\lambda=\omega_N\circ\phi$). Put
$
H_1=i^\vee\circ\lambda\circ i:
J^{p\text{-new}}\longrightarrow(J^{p\text{-new}})^\vee.
$

\begin{lem}\label{lem:pnew-selfdual}
	If $\mathfrak n$ is not $p$-old, the $\mathfrak m$-primary part of
	$\ker H_1$ vanishes. Consequently,
	$
	H_1^*:\mathcal X_p((J^{p\text{-new}})^\vee)_{\mathfrak m}
	\xrightarrow{\ \sim\ }
	\mathcal X_p(J^{p\text{-new}})_{\mathfrak m}
	$
	is a $\mathbb T^{p\text{-new}}(N)_{\mathfrak m}$-equivariant isomorphism.
\end{lem}

\begin{proof}
	The dual of the inclusion $i$ is the optimal quotient
	$
	i^\vee:J^\vee\twoheadrightarrow(J^{p\text{-new}})^\vee
	\simeq J^\vee/IJ^\vee.
	$
	Since $\lambda$ is Hecke-equivariant,
	$\lambda^{-1}(IJ^\vee)=IJ$. Thus
	$
	\ker H_1=J^{p\text{-new}}\cap IJ.
	$
	$IJ$ is the $p$-old subvariety and $J^{p\text{-new}}\cap IJ$ is finite. Its Hecke support lies in the intersection of
	the $p$-old and $p$-new supports. It therefore has no
	$\mathfrak m$-primary part when $\mathfrak n$ is not $p$-old.
	The induced map on localized character groups is then an isomorphism.
\end{proof}

\begin{proof}[Proof of Theorem~\ref{thm:pnew-component}]
	The injection of the $p$-new abelian subvariety induces an injection on
	the tori of special fibers, hence a surjection
	$
	i^*:\mathcal X_p(J)\twoheadrightarrow
	\mathcal X_p(J^{p\text{-new}})
	$ (see
	\cite[Section~5]{helm2007maps}). By
	\cite[Theorem~3.10]{ribet1990modular}, the Hecke action on
	$\mathcal X_p(J)$ factors through $\mathbb T^{p\text{-new}}(N)$.
	Moreover, \cite[Corollary~5.3]{helm2007maps} identifies the character
	group of $J^{p\text{-new}}$ with the quotient of
	$\mathcal X_p(J)$ by $T^{p\text{-new}}(N)$, modulo $\mathbb Z$-torsion. Consequently these two
	character groups have the same $\mathbb Z$-rank. Since they are free
	abelian groups, the surjection $i^*$ is an isomorphism.
	Now suppose $\mathfrak n$ is not $p$-old. By
	Lemma~\ref{lem:pnew-selfdual}, $H_1^*$ is an isomorphism on localized
	character groups, so $J^{p\text{-new}}$ is locally self-dual at $\mathfrak m$.
	Since $J/J^{p\text{-new}}$ is $p$-old, there exists
	$t\in\mathbb T(N)\setminus\mathfrak n$ which annihilates
	$J/J^{p\text{-new}}$. Thus the action $t_J:J\to J$ defined by $t$, factors through
	$i$, giving a homomorphism $h:J\to J^{p\text{-new}}$ such that
	$
	i\circ h=t_J.
	$
	Since $i$ is Hecke-equivariant and injective, restricting to
	$J^{p\text{-new}}$ also gives
	$
	h\circ i=t_{J^{p\text{-new}}}.
	$
	Dualizing, we obtain
	$
	h^\vee\circ i^\vee=t_{J^\vee},~
	i^\vee\circ h^\vee=t_{(J^{p\text{-new}})^\vee}.
	$
	As $t\notin\mathfrak n$, its action becomes invertible after
	localization at $\mathfrak n$ (and at $\mathfrak m$ on the $p$-new
	quotient). The induced maps on component groups therefore show that
	$
	(i^\vee)_*:\Phi_p(J^\vee)_{\mathfrak n}
	\xrightarrow{\ \sim\ }
	\Phi_p((J^{p\text{-new}})^\vee)_{\mathfrak m},
	$
	as required.
\end{proof}

		\section{Component-group surjectivity and level lowering}
\label{sec:level-lowering}

We now study the relation between component-group surjectivity and level
lowering. Write
$
B=J^{p\text{-new}}=J_0(N)^{p-\text{new}}, ~ A=A_f\subseteq B$,
where $f$ is a $p$-new normalized eigenform of level $N$.
Here $A$ denotes the Hecke subvariety dual to the corresponding
optimal quotient. We use the optimal
dual quotient
$
\pi:B^\vee\twoheadrightarrow A^\vee.
$
Let $\mathbb T_A$ be the image of $\mathbb T(N)$ in
$\operatorname{End}_{\mathbb Q}(A)$, let $\mathfrak m$ be a maximal
ideal of $\mathbb T_A$ of residue characteristic $\ell$, and write
$\mathfrak n$ for its preimage in $\mathbb T^{p\text{-new}}(N)$. We denote its
preimage in $\mathbb T^{p\text{-new}}(N)$ by the same symbol
$\mathfrak m$, where this causes no ambiguity.

\begin{thm}[Component-group criterion for level lowering]
	\label{thm:component-group-level-lowering}
	Suppose $\bar\rho_{\mathfrak n}$ is absolutely irreducible and $\bar\rho_{\mathfrak n}$ is
	finite at $p$. If
	$
	\pi_*:\Phi_p(B^\vee)_{\mathfrak m}
	\twoheadrightarrow\Phi_p(A^\vee)_{\mathfrak m}
	$
	is surjective, then $\mathfrak n$ is $p$-old.
\end{thm}

\begin{proof}
	Suppose that $\mathfrak n$ is not $p$-old. By
	Theorem~\ref{thm:pnew-component},
	$
	\Phi_p(J_0(N)^\vee)_{\mathfrak n}
	\xrightarrow{\sim}\Phi_p(B^\vee)_{\mathfrak m}.
	$
	The group $\Phi_p(J_0(N))$ is Eisenstein by
	\cite[Theorem~3.12]{ribet1990modular}. Absolute irreducibility of
	$\bar\rho_{\mathfrak n}$ therefore implies
	$\Phi_p(J_0(N))_{\mathfrak n}=0$. The $\mathbb{T}(N)$-equivariant isomorphism $\lambda: J \rightarrow J^\vee$ of Section 4, identifies the same localization of the
	component group of $J_0(N)^\vee$ with zero. Hence
	$\Phi_p(B^\vee)_{\mathfrak m}=0$, and surjectivity gives
	$
	\Phi_p(A^\vee)_{\mathfrak m}=0.
	$
	This contradicts the fact that $\Phi_p(A^\vee)_{\mathfrak m}=0$ if and only if  	$\bar\rho_{\mathfrak n}$ is not finite flat (see \cite{khare2003isomorphisms}, Proposition 3, Part 1, for example). Thus $\mathfrak n$ must be $p$-old.
\end{proof}

In section 5.2, we prove the local surjectivity in the above theorem under the DVR condition \ref{eq:coefficient-DVR} by using results of Helm in \cite{helm2007maps}, which itself uses some crucial techniques in Ribet's level lowering results~(Corollary 5.10). This shows that the level lowering and the surjectivity problem are (at least)in many cases equivalent.
\subsection{Mazur's Principle}
In this section, using the results of previous sections, we formulate an alternative proof of Mazur's principle away from the case $\ell^2 \mid N$. 
With notations as in the previous section, we begin with recalling the Mazur's principle (\cite[Theorem~6.1]{ribet1990modular} and \cite[Theorem~2.8]{edixhoven1992weight}).

\begin{thm}[Mazur's principle]
	\label{thm:component-group-level-lowering}
	Suppose $\bar\rho_{\mathfrak n}$ is absolutely irreducible, and 
	$\bar\rho_{\mathfrak n}$ is finite at $p$. If $\ell\nmid p-1$
	then $\mathfrak n$ is $p$-old.
\end{thm} 
We use the criterion 5.1 to give an alternative proof of this theorem and recover it except for the $\ell^2 \mid N$ case. First we use the following result of Emerton in \cite{emerton2003optimal}~(cf. \cite[Theorem D]{Hamidi2026Comparison}). 

\begin{prop}
	
	If \label{prop:one-sided-pnew}
	$
	\mathcal X_p(B)_{\mathfrak m}
	\cong \mathbb T^{p\text{-new}}(N)_{\mathfrak m}
	$
	as modules over the localized $p$-new Hecke algebra, then
	$
	\Phi_p(B^\vee)_{\mathfrak m}
	\longrightarrow\Phi_p(A^\vee)_{\mathfrak m}
	$
	is surjective.
\end{prop}

\begin{proof}
	See Corollary 3.6 of \cite{emerton2003optimal}.
\end{proof}
Write $N=pM$, with $p\nmid M$. Theorem~\ref{thm:pnew-component}
identifies the character group of $B$ with
$\mathcal X_p(J_0(pM))$. The following lemma provides the multiplicity one we need.

\begin{lem}\label{lem:character-multiplicity-one}
	Suppose that $\mathfrak n$ is non-Eisenstein and not $p$-old.
	If $\ell\nmid p-1$ and $\ell\nmid N$, then
	$
	\mathcal X_p(J)_{\mathfrak n}\cong\mathbb T^{p\text{-new}}(N)_{\mathfrak n}.
	$
\end{lem}
\begin{proof}
	Since $p\mid N$, the hypotheses imply $\ell\nmid 2p(p-1)N$.
	By \cite[Theorem~6.4]{ribet1990modular},
	$
	\dim_{\mathbb T/\mathfrak n}
	\mathcal X_p(J)/\mathfrak n\mathcal X_p(J)\leq 1.
	$
	The result follows from Nakayama's lemma and the faithfulness of the $p$-new Hecke action, since $\mathfrak n$ is not $p$-old.
\end{proof}
We put everything together to obtain the following corollary.
\begin{cor}[Mazur's principle away from $\ell^2\mid N$]\label{cor:mazur-away-square}
	Suppose that $\bar\rho_{\mathfrak n}$ is absolutely irreducible and finite at $p$. If $\ell\nmid p-1$ and $\ell^2\nmid N$, then $\mathfrak n$ is $p$-old.
\end{cor}
\begin{proof}
	First assume $\ell\nmid N$. Suppose that $\mathfrak n$ is not $p$-old. By Theorem~4.1,
	$
	\Phi_p(J^\vee)_{\mathfrak n}\cong\Phi_p(B^\vee)_{\mathfrak m}.
	$
	On the other hand, Lemma~\ref{lem:character-multiplicity-one} and Proposition~5.3 imply that
	$
	\Phi_p(B^\vee)_{\mathfrak m}\longrightarrow\Phi_p(A^\vee)_{\mathfrak m}
	$
	is surjective. Since $\bar\rho_{\mathfrak n}$ is finite at $p$, we have $\Phi_p(A^\vee)_{\mathfrak m}\neq 0$, whereas
	$
	\Phi_p(J^\vee)_{\mathfrak n}\cong\Phi_p(B^\vee)_{\mathfrak m}=0,
	$
	because $\mathfrak n$ is non-Eisenstein and $\Phi_p(J^\vee)$ is Eisenstein. This is a contradiction.
	Now suppose $\ell\mid N$. Since $\ell^2\nmid N$, there are two cases.
	If $\mathfrak n$ is $\ell$-old, then either $\ell=p$, in which case the conclusion is immediate, or $\ell\neq p$. In the latter case, we work at level $N/\ell$, which is prime to $\ell$. Applying the case $\ell\nmid N$ at this lower level shows that $\mathfrak n$ is $p$-old.
	If $\mathfrak n$ is not $\ell$-old, then $\bar\rho_{\mathfrak n}$ is not modular at level $N/\ell$. By the main theorem of \cite{mazur1991two},
	$
	\dim_{\mathbb T/\mathfrak n}J[\mathfrak n]=2.
	$
	Using the Grothendieck's monodromy (\cite[Theorem~6.4]{ribet1990modular}), we obtain
	$
	\dim_{\mathbb T/\mathfrak n}
	\mathcal X_p(J)/\mathfrak n\mathcal X_p(J)\leq 1.
	$
	For $\ell=p$, including $\ell=p=2$, the corresponding character-group bound follows from the monodromy argument of \cite[Propositions~17 and~18]{mazur1991two}. Hence, by Nakayama's lemma,
	$
	\mathcal X_p(J)_{\mathfrak n}\cong\mathbb T(N)_{\mathfrak n}.
	$
	The same argument as in the case $\ell\nmid N$ completes the proof.
\end{proof}

\subsection{Component-Group Surjectivity: A Domino Principle}

In this section, we establish a surjectivity result for the maps on
component groups induced by optimal quotients of new subvarieties of
Jacobians of Shimura curves, localized at non-Eisenstein maximal ideals,
under an additional rank one hypothesis, which is always satisfied for the case of elliptic curves.

Let $N=DM$ be a square-free positive integer, where $D$ is the product
of an even number of primes. Let $J^D(M)$ denote the Jacobian of the
Shimura curve over $\mathbb{Q}$ of discriminant $D$ and level $M$.
Fix a divisor $Q$ of $M$ and a prime $p$ dividing $QD$.
Let $\mathbb{T}^D_Q(M)$ denote the $Q$-new quotient of the Hecke
algebra $\mathbb{T}^D(M)$ acting on the $Q$-new subvariety $J^D_Q(M)$, and let
$\mathfrak{m}$ be a maximal ideal of $\mathbb{T}^D_Q(M)$
of residue characteristic $\ell$.
Let $f$ be a newform of weight $2$ and level $N$, and let
$I=I_f$ denote the augmentation ideal associated with $f$ in
$\mathbb{T}^D_Q(M)$. 
Let $
A^D(M)^\vee=J^D(M)/IJ^D(M)
$
be the corresponding optimal quotient of $J^D(M)$, where
$A^D(M)\subseteq J^D(M)$ is the dual abelian subvariety.

Denote the quotient map by
$
\pi:J^D(M)^\vee\longrightarrow A^D(M)^\vee,
$
and let
$
\pi_*:\Phi_p(J^D(M)^\vee)\longrightarrow\Phi_p(A^D(M)^\vee)
$
be the induced map on component groups. Ribet-Takahashi, used the level lowering techniques of Ribet in \cite{ribet1990modular} to prove the surjectivity of $\pi_*$ locally at non-Eisenstein primes when $D>1$, $A^D(M)$ is an elliptic curve, and $M$ is not prime.

\begin{prop}[Ribet-Takahashi]
	Suppose $D>1$, $A^D(M)$ is an elliptic curve, and $M$ is not  prime. If $\mathfrak{m}$ is non-Eisenstein, then $	\operatorname{coker}(\pi_*)_{\mathfrak{m}}=0$.
\end{prop}

The above result was generalized by Khare
\cite{khare2003isomorphisms} to higher-dimensional modular abelian
varieties under additional hypotheses, including a local rank-one
condition and restrictions on the residual representation.
For elliptic curves, Böckle--Khare--Manning
\cite{bockle2024wiles} removed the restrictions that $M$ be composite
and that the conductor be squarefree, under a stronger irreducibility
hypothesis and assuming that the residual characteristic is odd
and prime to the conductor (see Corollaries 7.9 and 7.10 of
$loc.~cit.$).

Here we remove the $M$ being prime condition of Ribet-Takahashi for $\ell>3$ in Helm's auxiliary-level setting in \cite{helm2007maps} and also for $\ell=3$ in the squarefree $\Gamma_0$-level setting, and extend these results to new quotients, including the case $D=1$.
In each of the situations below, let $J$ denote the indicated Jacobian or
new quotient, and suppose that $\pi$ factors through the canonical map
$J^\vee\longrightarrow A^D(M)^\vee$. We continue to write $\pi_*$ for the
induced map on component groups
$
\pi_*:\Phi_p(J^\vee)\longrightarrow\Phi_p(A^D(M)^\vee).
$
The main result of this section is the following which exhibits a domino phenomenon: a suitable condition at a single prime \(q\) implies the surjectivity of component-group maps at every relevant prime \(p\).

We use an arithmetic variant of the localized endomorphism comparison
from Theorem~\ref{thm:rigid-endomorphism-comparison}. Unlike that
theorem, the following consequence of Helm's work does not require
local self-duality or a hypothesis on roots of unity. It is specific
to the Hecke-isogeny categories considered in
\cite{helm2007maps} and requires controllability and the exclusion
of Helm's exceptional maximal ideals, with the extension to odd
residue characteristic in \cite{Hamidi2026Comparison}.

\begin{lem}[Endomorphism comparison under controllability]
	\label{lem:helm-controllable-endomorphisms}
	Let $A$ be an abelian variety that belongs to one of Helm's Hecke-isogeny categories, and
	suppose that $A$ has purely toric reduction at $r$. Let
	$\mathfrak m$ be non-Eisenstein and controllable at $r$ in the sense of Definition 6.4 of \cite{helm2007maps}. Then, with the
	canonical dual Hecke action on $A^\vee$, the natural map is an isomorphism:
	\[
	\operatorname{End}_{\mathbb T}(A^\vee)_{\mathfrak m}
	\xrightarrow{\ \sim\ }
	\operatorname{End}_{\mathbb T_{\mathfrak m}}
	(\mathcal X_r(A)_{\mathfrak m}).
	\]
\end{lem}

\begin{proof}
	Write $R=\mathbb T_{\mathfrak m}$ and
	$[A]=\operatorname{Hom}_{\mathbb T}(A,J_{\min})_{\mathfrak m}$,
	where $J_{\min}$ is Helm's minimal object. By
	\cite[Theorem~4.11 and Proposition~4.14]{helm2007maps},
	$
	\operatorname{End}_{\mathbb T}(A)_{\mathfrak m}^{\mathrm{op}}
	\cong\operatorname{End}_{R}([A]),
	~
	\mathcal X_r(A)_{\mathfrak m}
	\cong\mathcal X_r(J_{\min})_{\mathfrak m}\otimes_R[A].
	$
	Controllability implies
	$\mathcal X_r(J_{\min})_{\mathfrak m}\cong R$ by
	\cite[Lemma~6.5]{helm2007maps}. Thus the natural action
	identifies $\operatorname{End}_{\mathbb T}(A)_{\mathfrak m}^{\mathrm{op}}$
	with $\operatorname{End}_{R}(\mathcal X_r(A)_{\mathfrak m})$.
	Dualizing algebraic endomorphisms yields the assertion.
\end{proof}

In the rest of the section we obtain an application to the local component-group-map surjectivity.

\paragraph{Controllability.}
We recall the local condition of \cite[Definition~6.4]{helm2007maps}.
Let $\mathfrak m$ be a non-Eisenstein maximal ideal of a Hecke algebra
attached to one of the relevant Hecke-isogeny categories,
let $\ell$ be its residue characteristic, and write
$\overline\rho_{\mathfrak m}$ for the associated absolutely irreducible
residual representation. For a prime $q$ in the squarefree level,
$\mathfrak m$ is said to be \emph{controllable at $q$} if either of the
following conditions holds:
\begin{enumerate}
	\item $\overline\rho_{\mathfrak m}$ is not finite at $q$;
	\item $q\ne\ell$, $\overline\rho_{\mathfrak m}$ is unramified at
	$q$, and $\overline\rho_{\mathfrak m}(\operatorname{Frob}_q)$ is
	not scalar;
	\item $q=\ell\ne2$;
	\item $q=\ell=2$ and the image of a decomposition group at $2$
	under $\overline\rho_{\mathfrak m}$ is not contained in the
	scalar matrices.
\end{enumerate}
Here ``finite at $q$'' means that the Galois representation $\overline\rho_{\mathfrak m}$
extends to a finite flat group scheme over $\mathbb Z_q$.

\begin{rem}[A convenient numerical criterion]\label{rem}
	The condition $\ell\nmid(q-1)$ implies controllability at $q$
	for odd $\ell$. Indeed, for $q\ne\ell$, failure of controllability
	would imply that the residual representation is unramified with scalar
	Frobenius at $q$; in this setting the scalar is given by the local
	Atkin--Lehner operator and has square $1$. Comparing determinants with
	the mod-$\ell$ cyclotomic character then gives
	$q\equiv1\pmod\ell$. If $q=\ell$ is odd, controllability holds by
	the third clause of the definition.
	The controllability hypothesis is therefore often satisfied and is
	strictly more flexible than this numerical criterion. It also makes
	sense for $\ell=2$, whereas always $2\mid(q-1)$ when $q$ is odd.
\end{rem} 

\begin{thm}[Domino Principle]\label{domino}
	Fix a prime $p$. Let $J$ be one of the abelian varieties specified
	below, with its Hecke algebra $\mathbb T$, and let $I$ be the ideal
	defining the optimal quotient $A^D(M)$. Let $\mathfrak m$ be a
	non-Eisenstein maximal ideal of $\mathbb T$ containing $I$, with
	residue characteristic $\ell$. Assume that the completion of
	$(\mathbb T/I)_{\mathfrak m}$ is a discrete valuation ring
	(for instance, when $A^D(M)$ is an elliptic curve).
	Then
	$
	\operatorname{coker}(\pi_*)_{\mathfrak m}=0
	$
	in each of the following cases:
	\begin{enumerate}
		\item $J=J_Q^D(M)$, $\mathbb T=\mathbb T_Q^D(M)$, and
		$p\mid D$, provided either $\ell>2$ or the inverse image of
		$\mathfrak m$ in the Hecke algebra of $J^D(M)$ is
		controllable at some prime $r\mid D$;
		
		\item $J=J_Q^D(M)$, $\mathbb T=\mathbb T_Q^D(M)$, and
		$p\mid Q$, provided that there exists a prime $q\mid Q$
		at which $\mathfrak m$ is controllable;
		
		\item $J=J_p^D(M)$, $\mathbb T=\mathbb T_p^D(M)$, and
		$p\mid M$, provided that there exists a prime $q\mid M$ such that
		$\mathfrak m$ is not $q$-old. Here $q$-oldness is understood with respect to
		the full modular Jacobian $J_0(DM)$. 
	\end{enumerate}
	For $\ell>3$, all three cases also hold with auxiliary level
	structure $\Gamma_0(M)\cap\Gamma'$ as in \cite{helm2007maps},
	where the level of $\Gamma'$ is prime to $DM$, for the
	corresponding optimal quotients.
\end{thm}

\begin{proof}
	We begin with case~(1) for $Q=1$. For $\ell>3$, the endomorphism
	comparison of \cite[Proposition~8.13]{helm2007maps} applies also
	with auxiliary level structures. For $\ell=3$, the corresponding
	comparison for $J^D(M)$ follows from \cite{Hamidi2026Comparison}.
	In either case, taking the two varieties equal gives a natural isomorphism
	$
	\operatorname{End}_{\mathbb T}(J^D(M)^\vee)_{\mathfrak m}
	\xrightarrow{\ \sim\ }
	\operatorname{End}_{\mathbb T_{\mathfrak m}}
	\bigl(\mathcal X_p(J^D(M))_{\mathfrak m}\bigr).
	$
	If instead $\mathfrak m$ is controllable at some prime $r\mid D$,
	the localized assertion of Helm's Theorem~6.1 follows from
	\cite[Corollary~6.6]{helm2007maps}. The proof of
	\cite[Proposition~8.13]{helm2007maps} then applies verbatim:
	by \cite[Proposition~4.14]{helm2007maps}, the character group
	at $p$ is obtained by tensoring the character group of the
	minimal object with its associated Hecke module, and the
	Hom--tensor adjunction, together with
	\cite[Theorem~4.11]{helm2007maps}, gives the
	above endomorphism comparison. For $\ell=2$, the needed use of
	Proposition~5.7 in Corollary~6.6 remains valid for the $\Gamma_0$-type levels; see
	\cite[Remark~8.12]{helm2007maps}. In either situation,
	Proposition~2.3 gives
	$\operatorname{coker}(\pi_*)_{\mathfrak m}=0$.
	For general $Q$, the natural inclusion of the $Q$-new subvariety
	into $J^D(M)$ induces a quotient map on duals
	$
	\alpha:J^D(M)^\vee\longrightarrow J_Q^D(M)^\vee.
	$
	Since $A^D(M)^\vee$ is a quotient of the $Q$-new part, the quotient
	map from $J^D(M)^\vee$ to $A^D(M)^\vee$ factors through $\alpha$.
	Let $\mathfrak m_1$ be the inverse image of $\mathfrak m$ in
	the Hecke algebra of $J^D(M)$. The coefficient order of
	$A^D(M)$, hence the discrete-valuation-ring hypothesis, is
	unchanged. Thus
	\[
	\Phi_p(J^D(M)^\vee)_{\mathfrak m_1}
	\xrightarrow{\alpha_*}
	\Phi_p(J_Q^D(M)^\vee)_{\mathfrak m}
	\xrightarrow{\pi_*}
	\Phi_p(A^D(M)^\vee)_{\mathfrak m}.
	\]
	The composite is surjective by the case $Q=1$, so $\pi_*$ is
	surjective. This proves case~(1).
	
	We turn to case~(2). If $q=p$, controllability at $p$
	and Lemma~\ref{lem:helm-controllable-endomorphisms} give
	$
	\operatorname{End}_{\mathbb T}(J^\vee)_{\mathfrak m}
	\xrightarrow{\ \sim\ }
	\operatorname{End}_{\mathbb T_{\mathfrak m}}
	\bigl(\mathcal X_p(J)_{\mathfrak m}\bigr).
	$
	Proposition~2.3 again proves the assertion. Henceforth we may
	take $q\ne p$.
	We first prove case~(2) for $Q=pq$. Set
	$
	U=J_{pq}^D(M),~
	V=J^{Dpq}(M/pq),~
	\mathbb T=\mathbb T_{pq}^D(M).
	$
	Via the Jacquet--Langlands correspondence, the same Hecke algebra
	acts on $U$ and $V$. Since $\mathfrak m$ is non-Eisenstein,
	\cite[Proposition~5.7]{helm2007maps}, applied in both orders,
	yields
	$
	\mathcal X_p(U)_{\mathfrak m}
	\cong\mathcal X_q(V)_{\mathfrak m}^{*},
	~
	\mathcal X_q(U)_{\mathfrak m}
	\cong\mathcal X_p(V)_{\mathfrak m}^{*}.
	$
	For $r\in\{p,q\}$ and $W\in\{U,V\}$, write
	$
	E_r(W)=\operatorname{End}_{\mathbb T_{\mathfrak m}}
	\bigl(\mathcal X_r(W)_{\mathfrak m}\bigr),
	~
	F(W)=\operatorname{End}_{\mathbb T}(W^\vee)_{\mathfrak m}.
	$
	These character groups have generic rank one over the
	corresponding commutative Hecke algebra, so their equivariant
	endomorphism rings are commutative; $\mathbb Z_\ell$-duality
	therefore identifies the endomorphism rings of a character
	group and its dual.
	By controllability at $q$ and
	Lemma~\ref{lem:helm-controllable-endomorphisms},
	$
	E_q(U)\cong F(U),~ E_q(V)\cong F(V).
	$
	The isomorphisms above, together with the faithful action of
	algebraic endomorphisms on character groups at $p$, give
	\[
	\begin{aligned}
		E_p(U)
		&\xrightarrow{\ \sim\ } E_q(V)
		\xrightarrow{\ \sim\ } F(V) \\
		&\hookrightarrow E_p(V)
		\xrightarrow{\ \sim\ } E_q(U)
		\xrightarrow{\ \sim\ } F(U).
	\end{aligned}
	\]
	The rational rank-one property identifies all these
	endomorphism rings rationally with
	$\mathbb T_{\mathfrak m}\otimes\mathbb Q_\ell$.
	Since the maps in the above chain respect the Hecke
	action and restrict to the identity on this rational algebra,
	the resulting embedding $E_p(U)\hookrightarrow F(U)$ is
	inverse to the natural injection $F(U)\hookrightarrow E_p(U)$.
	Consequently every Hecke-equivariant endomorphism of
	$\mathcal X_p(U)_{\mathfrak m}$ is algebraic after localization.
	Proposition~2.3 therefore yields
	$\operatorname{coker}(\pi_*)_{\mathfrak m}=0$ for $Q=pq$.
	For arbitrary $Q$ with $pq\mid Q$, use the quotient on duals
	$
	\alpha:J_{pq}^D(M)^\vee\longrightarrow J_Q^D(M)^\vee.
	$
	If $\mathfrak m_{pq}$ denotes the inverse image of
	$\mathfrak m$ in $\mathbb T_{pq}^D(M)$, the component-group
	maps factor as
	\[
	\Phi_p(J_{pq}^D(M)^\vee)_{\mathfrak m_{pq}}
	\xrightarrow{\alpha_*}
	\Phi_p(J_Q^D(M)^\vee)_{\mathfrak m}
	\xrightarrow{\pi_*}
	\Phi_p(A^D(M)^\vee)_{\mathfrak m}.
	\]
	The composite is surjective by the $Q=pq$ case, so the second
	arrow is surjective. For $\ell>3$, the same argument applies
	with auxiliary level $\Gamma'$, using the corresponding
	character-group comparisons of \cite{helm2007maps}.
	This proves case~(2).

	Finally, consider case~(3). Since $\mathfrak m$ is not
	$q$-old, Serre's modularity conjecture, proved by
	Khare--Wintenberger \cite{kharewintenberger2009serre}, implies that
	$\bar\rho_{\mathfrak m}$ is not finite at $q$.
	In particular, $\mathfrak m$ is controllable at $q$.
	If $q=p$, the assertion follows from case~(2).
	Henceforth assume $q\ne p$, and set
	$
	C=J_p^D(M),~ U=J_{pq}^D(M).
	$
	Let $i:U\hookrightarrow C$ be the natural inclusion of the
	$q$-new subvariety. Because $\mathfrak m$ is not $q$-old,
	the $q$-old contribution vanishes after localization at
	$\mathfrak m$. The argument used in the proof of Theorem \ref{thm:pnew-component} gives an isomorphism
	$
	(i^\vee)_*:
	\Phi_p(C^\vee)_{\mathfrak m}
	\xrightarrow{\ \sim\ }
	\Phi_p(U^\vee)_{\mathfrak m}.
	$
	Here the maximal ideal on the right is its image in the
	$pq$-new Hecke algebra. The optimal quotient map factors
	through $U^\vee$ after localization, giving
	$
	\Phi_p(C^\vee)_{\mathfrak m}
	\xrightarrow[\sim]{\ (i^\vee)_*\ }
	\Phi_p(U^\vee)_{\mathfrak m}
	\longrightarrow
	\Phi_p(A^D(M)^\vee)_{\mathfrak m}.
	$
	Since $\mathfrak{m}$ is controllable at $q$, invoking the case~(2), we conclude that the second map above is surjective while the first
	arrow is an isomorphism. For $\ell>3$, the argument also works
	with auxiliary level $\Gamma'$. Indeed, let $C$ and $U$ denote
	the corresponding $p$-new and $pq$-new subvarieties at that  level. Since $U$ is the $q$-new subvariety of $C$,
    the quotient $C/U$ is Hecke-equivariantly isogenous
    to the $q$-old part of $C$.
    Since $\mathfrak m$
	is not $q$-old, there exists $t\in\mathbb T\setminus\mathfrak m$
	which annihilates $C/U$. Consequently, the endomorphism
	$t_C:C\to C$ (the restriction map) factors through $i$, giving a homomorphism
	$h:C\to U$ with
	$
	i\circ h=t_C,~ h\circ i=t_U.
	$
	Dualizing yields
	$
	h^\vee\circ i^\vee=t_{C^\vee},~
	i^\vee\circ h^\vee=t_{U^\vee}.
	$
	As $t$ is invertible after localization at $\mathfrak m$,
	these identities show that
	$
	(i^\vee)_*:\Phi_p(C^\vee)_{\mathfrak m}
	\xrightarrow{\ \sim\ }\Phi_p(U^\vee)_{\mathfrak m}
	$
	is an isomorphism. The quotient map to $A^D(M)$ factors
	through $U^\vee$ after localization, and case~(2), already
	established at auxiliary level for $\ell>3$, makes the
	resulting component-group map surjective. Hence
	$\operatorname{coker}(\pi_*)_{\mathfrak m}=0$.
\end{proof}

From part 3 of the above theorem immediately follows that:
\begin{cor}
	Let $M$ be squarefree, let $p\mid M$, and retain the
	hypotheses of Theorem~\ref{domino}. Then
	the map
	$
	\pi_*:\Phi_p((J^1_p(M))^\vee)_{\mathfrak m}
	\longrightarrow
	\Phi_p(A^1(M)^\vee)_{\mathfrak m}
	$
	is surjective.
\end{cor}

\begin{proof}

There exists a prime $q\mid M$ such that
$\mathfrak m$ is not $q$-old. Indeed, otherwise,
by Serre's modularity conjecture
\cite{kharewintenberger2009serre},
$\bar\rho_{\mathfrak m}$ would arise from a
weight-two modular form of level $1$, which is
impossible. Hence, the assertion follows
from case~(3).

\end{proof}

In Section 9 of \cite{helm2007maps}, Helm conjectures that Theorem 6.1 of $loc.~cit.$
remains valid at all non-Eisenstein maximal ideals for certain
congruence-level structures. The odd residue characteristic case is
established in \cite{Hamidi2026Comparison} in our $\Gamma_0$-setting. Since Helm's Theorem 6.1 implies Proposition~8.13 of $loc.~cit.$, which is
used in our preceding result, his conjecture provides motivation for
expecting surjectivity also in the remaining exceptional characteristic
$2$.
Regarding the DVR hypothesis, an examination of the data in the LMFDB \cite{lmfdb}
shows that, among weight-two newforms with trivial nebentypus and level
at most \(1000\), only relatively few have coefficient rings that are
nonmaximal at odd primes. Nonmaximality occurs more frequently at $2$.
Nevertheless, we expect the surjectivity phenomenon to persist even
when the localized coefficient order is not a discrete valuation ring. Therefore, we propose the following conjecture.

\begin{conj}\label{conj}
	Let $J$ be the $Q$-new subvariety of the Jacobian of a
	Shimura curve of discriminant $D$ and level
	$\Gamma_0(M)\cap\Gamma'$, where $M$ is squarefree,
	$Q\mid M$, and $\Gamma'$ is an auxiliary congruence subgroup
	as in \cite{helm2007maps}. Let $p$ be a prime dividing $QD$.
	Let $\mathbb T$ be the Hecke algebra acting on $J$, and let
	$A\subseteq J$ be a Hecke subvariety defined by an ideal
	$I\subseteq\mathbb T$. Let
	$
	\pi:J^\vee\longrightarrow A^\vee
	$
	be the dual optimal quotient.
	If $\mathfrak m$ is a non-Eisenstein maximal ideal of
	$\mathbb T$ containing $I$, then
	$
	\operatorname{coker}(\pi_*)_{\mathfrak m}=0.
	$
\end{conj}

\begin{rem}
	The non-Eisenstein hypothesis in the preceding conjecture
	is necessary. Indeed, counterexamples to surjectivity
	at Eisenstein maximal ideals were found in
	\cite{Hamidi2026Comparison}.
\end{rem}

It is natural to ask whether the same surjectivity statement holds for the more general abelian varieties $J$ considered in Section 2.1, when $J$ is defined over $\mathbb Q$. Papikian and Rabinoff \cite{papikian2016optimal} construct, for every odd prime $p$, examples in which the map on component groups at $p$ fails to be surjective. Their examples, however, involve reducible mod-$2$ Galois representations. We modify their construction to obtain, for every odd prime $p$, infinitely many examples with an absolutely irreducible mod-$2$ representation. Thus the analogous statement for the general class of varieties in Section 2.1 fails even at non-Eisenstein maximal ideals. In particular, a nontrivial $2$-primary cokernel need not arise
from the existence of nontrivial rational $2$-torsion on the
optimal quotient.

\begin{exa}\label{exa}
	Fix an odd prime $p$. We construct a Jacobian $J/\mathbb Q$, two optimal quotients $\pi_i:J\to E_i$, and a non-Eisenstein maximal ideal $\mathfrak m$ for which
	$
	\Phi_p(J)=0,~
	\Phi_p(E_i)\cong\mathbb Z/2\mathbb Z,~\text{and}~
	\operatorname{coker}((\pi_i)_*)_{\mathfrak m}\ne0
	~(i=1,2).
	$
	The construction has six steps. Its starting points are the elliptic curves
	$
	y^2=x(x-1)(x-p),
	$
	and
	$
	y^2=x(x-1)(x+p)
	$
	from Section 5.3 of \cite{papikian2016optimal}. Gluing these curves along their $2$-torsion gives the desired behavior of component groups, but their $2$-torsion representations are reducible. We retain the gluing method while replacing the curves.
\end{exa}

\begin{proof}[Construction]
	\smallskip
	\noindent\textbf{Step 1: A curve with absolutely irreducible $2$-torsion.}
	Choose a nonzero \emph{odd} integer $s$ with $p\nmid s$, and set
	$
	f_s(x)=x(x-1)(x-p)+p^3s,
	~ E_1:\ y^2=f_s(x).
	$
	The polynomial $f_s$ has no rational root: since it is monic, any rational root would be an integer $r$, and then
	$
	r(r-1)(r-p)=-p^3s,
	$
	which is impossible because the left side is even and the right side is odd. Thus $f_s$ is irreducible over $\mathbb Q$.
	Put
	$
	D_s=(p-1)^2+2p(p+1)(2p-1)(p-2)s-27p^4s^2.
	$
	A direct calculation gives $\operatorname{disc}(f_s)=p^2D_s$. Since $p\ge3$ and $|s|\ge1$, we have $D_s<0$ since
	$
	(p-1)^2<p^4$ and
	$2p(p+1)(2p-1)(p-2)<4p^4
	$. In particular, the discriminant is not a square in $\mathbb Q$. The splitting field of the irreducible cubic $f_s$ therefore has Galois group $S_3$. As $S_3\cong\operatorname{GL}_2(\mathbb F_2)$ and its natural two-dimensional representation is absolutely irreducible, $E_1[2]$ is absolutely irreducible as a $G_{\mathbb Q}$-module.
	
	\smallskip
	\noindent\textbf{Step 2: Roots of $f_s$.}
	Write the roots of $f_s$ as $x_1,x_2,x_3$. We claim that they all belong to $\mathbb Q_p$ and may be labeled so that
	\[
	\begin{array}{c|ccc}
		&x_1&x_2&x_3\\ \hline
		x_i\pmod p&1&0&0\\
		v_p(x_i)&0&2&1
	\end{array}
	\qquad\text{and}\qquad x_1\equiv1\pmod{p^3}.
	\]
	Indeed, $f_s(x)\equiv x^2(x-1)\pmod p$, so Hensel's lemma gives a root $x_1\equiv1\pmod p$. From $x_1(x_1-1)(x_1-p)=-p^3s$, the two unit factors show that $v_p(x_1-1)=3$.
	For the other roots, put $x=pz$. The equation $f_s(x)=0$ becomes
	$
	F(z):=z(pz-1)(z-1)+ps=0,
	~ F(z)\equiv-z(z-1)\pmod p.
	$
	Both $0$ and $1$ are simple roots of $F$ modulo $p$. Hensel's lemma therefore gives roots $z_0\equiv0$ and $z_1\equiv1\pmod p$. The equation $F(z_0)=0$ shows that $v_p(z_0)=1$,
	while $z_1$ is a unit. Thus $x_2=pz_0$ and $x_3=pz_1$ have the claimed valuations. Denote the nonzero $2$-torsion points of $E_1$ by $P_i=(x_i,0)$.
	
	\smallskip
	\noindent\textbf{Step 3: Construct the second curve and the gluing map.}
	Let
	$
	g(x)=\frac{x(x-1)}p.
	$
	Define $E_2$ by the equation
	$y^2=\prod_{i=1}^3(x-q_i),$
	where $q_i=g(x_i)$.
	The valuations from Step 2 imply
	\[
	v_p(q_1)=2,\qquad v_p(q_2)=1,\qquad v_p(q_3)=0,
	\qquad q_3\equiv-1\pmod p.
	\]
	In particular, the $q_i$ are distinct and hence $E_2$ is an elliptic curve. To obtain an equation for $E_2$ over $\mathbb Q$, recall that
	$
	x_i(x_i-1)(x_i-p)=-p^3s.
	$
	Since $q_i=x_i(x_i-1)/p$, we obtain
	$
	q_i(x_i-p)=-p^2s,
	$
	and hence
	$
	x_i=p-\frac{p^2s}{q_i}.
	$ Substituting into $x_i(x_i-1)=pq_i$ gives
	$
	\left(p-\frac{p^2s}{q_i}\right)
	\left(p-1-\frac{p^2s}{q_i}\right)=pq_i.
	$
	Clearing denominators gives
	$
	q_i^3-(p-1)q_i^2+p(2p-1)sq_i-p^3s^2=0.
	$
	Consequently, $E_2$ is given by the equation
	$
	y^2=x^3-(p-1)x^2+p(2p-1)sx-p^3s^2.
	$
	Write $Q_i=(q_i,0)\in E_2[2]$. Since $g\in\mathbb Q(x)$, the assignment
	$
	\psi:E_1[2]\longrightarrow E_2[2],
	~ \psi(O)=O,~ \psi(P_i)=Q_i
	$
	is a $G_{\mathbb Q}$-equivariant isomorphism. Define its graph
	$
	G=\{(P,\psi(P)):P\in E_1[2]\}\subset E_1[2]\times E_2[2].
	$
	
	\smallskip
	\noindent\textbf{Step 4: Reduction and component groups.}
	The reductions of the equations of $E_1$ and $E_2$ modulo $p$ are
	$
	\widetilde E_1:\ y^2=x^2(x-1),
	$ and
	$\widetilde E_2:\ y^2=x^2(x+1).
	$
	Both have a node at $(0,0)$. Their standard Weierstrass invariants are
	\[
	\begin{array}{c|cc}
		&E_1&E_2\\ \hline
		c_4&16(p^2-p+1)&16((p-1)^2-3p(2p-1)s)\\[2pt]
		\Delta&16p^2D_s&16p^2s^2D_s
	\end{array}
	\]
	so $v_p(c_4(E_i))=0$ and $v_p(\Delta(E_i))=2$. Hence both curves have multiplicative reduction of type $I_2$, and
	$
	\Phi_p(E_1)\cong\Phi_p(E_2)\cong\mathbb Z/2\mathbb Z.
	$
	The reductions of $P_1$ and $Q_3$ lie in the smooth loci
	of the special fibers of $E_1$ and $E_2$, respectively,
	whereas the other nontrivial $2$-torsion points reduce
	to the nodes. Thus, on $2$-torsion, the kernels of the
	specialization maps to component groups are
	$\{O,P_1\}$ for $E_1$ and $\{O,Q_3\}$ for $E_2$
	(cf. Section 5.3 of \cite{papikian2016optimal}).

	\smallskip
	\noindent\textbf{Step 5: The glued surface is a Jacobian with trivial component group.}
	Set
	$J=(E_1\times E_2)/G$.
	Because $G$ is Galois-stable, $J$ is defined over $\mathbb Q$, and moreover, being isogenous to $E_1\times E_2$, it has purely toric reduction at $p$.
	By Proposition 5.2 and the discussion in Section 5.3 of
	\cite{papikian2016optimal}, to conclude that $J$ is a Jacobian, it remains to check that
	$E_1$ and $E_2$ are not geometrically isomorphic. For brevity, put
	\[
	A=(p-1)^2,\qquad B=3p(2p-1),\qquad C=p^2-p+1.
	\]
	The computations of $\Delta$ and $c_4$ above give
	\[
	j(E_1)=\frac{256C^3}{p^2D_s},
	\qquad
	j(E_2)=\frac{256(A-Bs)^3}{p^2s^2D_s},
	\qquad
	\frac{j(E_2)}{j(E_1)}=\frac{(A-Bs)^3}{s^2C^3}.
	\]
	Since $B-A-C=4p^2-2>0$ and $|s|\ge1$,
	\[
	|A-Bs|\ge B|s|-A\ge(B-A)|s|>C|s|.
	\]
	It follows that
	$
	\left|\frac{j(E_2)}{j(E_1)}\right|>|s|\ge1,
	$
	so $j(E_1)\ne j(E_2)$. Therefore $J$ is the Jacobian of a genus-two curve.
	It remains to compute its component group. Let
	$
	h:G\longrightarrow\Phi_p(E_1)\times\Phi_p(E_2)
	$
	be the canonical specialization map. Under the identification $G\cong E_1[2]$ induced by its first projection, the first component of $h$ has kernel $\{O,P_1\}$, whereas the second has kernel $\{O,P_3\}$, because $\psi(P_3)=Q_3$. Their intersection is trivial. Hence $h$ is injective. Since both sides have order four, it is an isomorphism. Lemma 5.3 of \cite{papikian2016optimal} now gives
	$
	\Phi_p(J)=0.
	$
	The natural quotient maps $\pi_i:J\to E_i$ have connected kernels, so they are optimal. Nonetheless, neither map
	$
	(\pi_i)_*:\Phi_p(J)\longrightarrow\Phi_p(E_i)
	$
	is surjective.
	
	\smallskip
	\noindent\textbf{Step 6: $\mathbb{T}$ and a non-Eisenstein ideal.}
	Consider the subring of $\operatorname{End}(J)$
	$
	\mathbb T=\{(a,b)\in\mathbb Z\times\mathbb Z:a\equiv b\pmod2\}.
	$
	The congruence condition ensures that the coordinatewise action on $E_1\times E_2$ preserves $G$ and hence descends to $J$. Let
	$
	\mathfrak m=\{(a,b)\in\mathbb T:a\equiv b\equiv0\pmod2\}.
	$
	$\mathfrak m$ is non-Eisenstein because its corresponding Galois representation is defined on $E_1[2]$ which is absolutely irreducible by Step 1.  
	Because $\Phi_p(J)=0$ while $\Phi_p(E_i)\cong\mathbb Z/2\mathbb Z$, and $\mathbb T$ acts on $E_i$ through its $i$th projection, we conclude that
	$
	\operatorname{coker}((\pi_i)_*)_{\mathfrak m}
	\cong\mathbb F_2\ne0
	~(i=1,2).
	$
	Finally, one has
	$
	\mathcal X_p(J)\otimes\mathbb Q
	\cong\mathbb T\otimes\mathbb Q
	$
	as $\mathbb T\otimes\mathbb Q$-modules. Thus $J$ belongs to the class considered in Section 2, as required. The formula for $j(E_1)$ is a nonconstant rational function of $s$. Since there are infinitely many odd integers $s$ prime to $p$, this yields infinitely many distinct examples for every odd prime $p$.
\end{proof}

\medskip
\noindent\textbf{Acknowledgements.}
I would like to express my sincere gratitude to Patrick Allen for
his guidance, support, and many valuable discussions throughout this
work. I am also grateful to Mihran Papikian for carefully reading an
earlier version of this manuscript and providing thoughtful comments.

		\bibliographystyle{plain}
		\bibliography{paper}

	\end{document}